\documentclass[english]{smfart}

\newtheorem{thm}{Theorem}[section]
\newtheorem{lemma}[thm]{Lemma}
\newtheorem{prop}[thm]{Proposition}
\newtheorem{cor}[thm]{Corollary}

\newtheorem{rmk}[thm]{Remark}
\newtheorem{example}[thm]{Example}

\newtheorem{defn}[thm]{Definition}

\numberwithin{equation}{thm}

\def \Pr {\textsl{Proof. - }}

\RequirePackage{amsmath}
\RequirePackage{amssymb}
\usepackage[english, francais]{babel}

\author[J. Sauloy]{Jacques Sauloy}
\address{Laboratoire Emile Picard, CNRS UMR 5580, U.F.R. M.I.G.,
118, route de Narbonne, 31062 Toulouse CEDEX 4}
\email{sauloy@picard.ups-tlse.fr}
\urladdr{http://picard.ups-tlse.fr/~sauloy/}
\title[Algebraic construction of the Stokes sheaf
for irregular $q$-difference
equations]{Algebraic construction of the  Stokes sheaf
for irregular linear \protect\boldmath$q$\protect\unboldmath-difference
equations}
\alttitle{Construction alg\'ebrique du faisceau de Stokes pour les
\'equations aux \protect\boldmath$
q$\protect\unboldmath-diff\'erences lin\'eaires
irr\'eguli\`eres}

\def\sq{\sigma_q}

\def\C{{\mathbf C}}
\def\Q{{\mathbf Q}}
\def\Z{{\mathbf Z}}
\def\R{{\mathbf R}}
\def\N{{\mathbf N}}

\def\F{{\mathcal{F}}}
\def\P{{\mathcal{P}}}
\def\B{{\mathcal{B}}}
\def\A{{\mathcal{A}}}
\def\O{{\mathcal{O}}}
\def\M{{\mathcal{M}}}
\def\Ma{{\text{Mat}}}

\def\D{{\mathcal{D}_{q,K}}}
\def\gr{{\text{gr}}}
\def\VV{{\mathcal{V}}}
\def\Eq{{\mathbf{E}_{q}}}
\def\EE{{\mathcal{E}}}
\def\G{{\mathfrak{G}}}
\def\g{{\mathfrak{g}}}
\def\U{{\mathfrak{U}}}
\def\V{{\mathfrak{V}}}
\def\k{{\C(z)}}
\def\Ra{{\C\{z\}}}
\def\Ka{{\C(\{z\})}}
\def\Rf{{\C[[z]]}}
\def\Kf{{\C((z))}}
\def\Rw{{\mathcal{O}(\C^{*})}}
\def\Kw{{\mathcal{M}(\C^{*})}}
\def\Rwg{{\mathcal{O}(\C^{*},0)}}
\def\Kwg{{\mathcal{M}(\C^{*},0)}}
\def\DM{{DiffMod(K,\sq)}}
\def\DMa{{DiffMod(\Ka,\sq)}}
\def\DMf{{DiffMod(\Kf,\sq)}}

\def\l{\left}
\def\r{\right}

\def\la{\lambda}
\def\La{\Lambda}
\def\th{\theta}
\def\Th{\Theta}

\def\div{\text{div}}

\begin{document}

\frontmatter

\begin{altabstract}
   La classification analytique locale
des \'equations aux $q$-diff\'erences lin\'eaires irr\'eguli\`eres
a \'et\'e r\'ecemment r\'ealis\'ee par J.-P. Ramis, J. Sauloy et
C. Zhang. Leur description fait intervenir un $q$-analogue du
faisceau de Stokes
et des th\'eor\`emes de type Malgrange-Sibuya et elle s'appuie
sur la sommation discr\`ete de C. Zhang. Nous montrons ici
comment retrouver une partie de ces r\'esultats par voie
alg\'ebrique et nous d\'ecrivons le 
d\'evissage $q$-Gevrey du $q$-faisceau de Stokes par des fibr\'es
vectoriels holomorphes sur une courbe elliptique.
\end{altabstract}

\begin{abstract}
   The local analytic classification of
irregular linear $q$-difference equations has recently been obtained
by J.-P. Ramis, J. Sauloy and C. Zhang. Their description involves
a $q$-analog of the Stokes sheaf and theorems of Malgrange-Sibuya type and is
based on a discrete summation process due to C. Zhang. We show here
another road to some of these results by algebraic means and we
describe the $q$-Gevrey devissage of the
$q$-Stokes sheaf by holomorphic vector bundles over an elliptic
curve.
\end{abstract}

\maketitle

\hfill{\it Je laisse aux nombreux avenirs (non \`a tous) }

\hfill{\it mon jardin aux sentiers qui bifurquent}

\hfill(Jorge Luis Borges, {\it Fictions}).

\tableofcontents

\mainmatter

%%%%%%%%%%%%%%%%%%%%%%%%%%%%%%%%%%%%%%%%%%%%%%%%%%%%%%%%%%%%%%%%%%%%%%%%%%%%%

% 1

\section{Introduction and general conventions}

% 1.1

\subsection{Introduction}

This paper deals with Birkhoff's program of 1941
(\cite{Birkhoff3}, see also \cite{Birkhoff1}) towards the local 
analytic classification
of $q$-difference equations and some extensions stated by J.-P.
Ramis in 1990 (\cite{RamisJPRTraum}).
\\

A full treatment of the Birkhoff program including the case of
irregular $q$-difference equations is being given in \cite{RSZ}.
The method used there closely follows the analytic
procedure developed in the last decades by B. Malgrange, Y.
Sibuya, J.-P. Ramis,\dots for the ``classical'' case, i.e., the
case of differential equations: adequate asymptotics, $q$-Stokes
phenomenon, $q$-Stokes sheaves and theorems of Malgrange-Sibuya
type; explicit cocycles are built using a \emph{discrete summation
process} due to C. Zhang (\cite{Zha02}) where the Jackson
$q$-integral and theta functions are introduced in place of the
Laplace integral and exponential kernels.

To get an idea of the classical theory for linear differential
equations one should look at the survey \cite{Varadarajan} by V.
S. Varadarajan, especially section 6, and to get some feeling of
how the change of landscape from differential equations to
$q$-difference equations operates, at the survey \cite{DRSZ} by
L. Di Vizio, J.-P. Ramis, J. Sauloy and C. Zhang.
\bigskip

The aim of this paper is to show how the harder analytic tools
can, to some extent, be replaced by much simpler algebraic
arguments. The problem under consideration being a transcendental
one we necessarily keep using analytic arguments but in their
most basic,  ``19th century style'', features only. In
particular, we avoid here using the discrete summation process.

Again our motivation is strongly pushed ahead by the
classical model of which we recall three main steps: the
{\it d\'evissage Gevrey} introduced by J.-P. Ramis
(\cite{RamisGevrey}) occured to be the fundamental tool for
understanding the Stokes phenomenon. The underlying algebra was
clarified by P. Deligne in \cite{Deligne86}, then put at work by
D. G. Babbitt and V. S. Varadarajan in \cite{BV} (see also
\cite{Varadarajan}) for moduli theoretic purposes. On the same
basis, effective methods, a natural summation and galoisian
properties were thoroughly explored by M. Loday-Richaud in
\cite{MLR}.

Some  specificities of our problem are due, on
one hand, to the fact that the sheaves to be considered are
quite similar to holomorphic vector bundles over an elliptic
curve, whence the benefit of GAGA theorems, and, on another hand,
  to the existence for
$q$-difference operators of an analytic factorisation without
equivalent for differential operators. Such a factorisation
originates in Birkhoff (\cite{Birkhoff3}), where it was
rather stated in terms of a triangular form
of the system. It has been revived by C. Zhang
(\cite{Zha99}, \cite{MZ}) in terms of factorisation and we will
use it in its linear guise, as a filtration of $q$-difference
modules (\cite{JSFIL}). \\

In this paper, following the classical theory recalled above, 
we build a $q$-Gevrey filtration on the
$q$-Stokes sheaf, thereby providing a $q$-analog of the Gevrey
devissage in the classical case. This $q$-devissage jointly with
a natural summation argument allows us  to prove  the $q$-analog
of a Malgrange-Sibuya theorem (theorem 3 of \cite{Varadarajan})
in  quite a direct and easy way; in particular, we avoid here
the Newlander-Nirenberg structural theorem used in \cite{RSZ}.
Our filtration is, in some way, easier to get than the classical
one : indeed, due to the forementionned canonical filtration of
$q$-difference  modules, our systems admit a natural
triangularisation which is independent of the choice of a Stokes
direction  and of the domination order of exponentials
(here replaced by theta functions). Also, our
filtration has a much nicer structure than the classical Gevrey
filtration since the so-called elementary sheaves of the
classical theory are here replaced by holomorphic vector bundles
endowed with a very simple structure over an elliptic curve
(they are tensor products of flat bundles by line bundles).

\bigskip

On the side of what this paper does not contain, there
is neither a study of confluency when $q$ goes to $1$,
nor any application to Galois theory. As for the former,
we hope to extend the results in \cite{JSAIF} to the
irregular case, but this seems a difficult matter. Only partial
results by C. Zhang are presently available, on significant
examples. As for the latter, it is easier to obtain as a
consequence of the present results that, under natural
restrictions, ``canonical Stokes operators are Galoisian'' like in
\cite{MLR}. However, to give this statement its full meaning, 
we have to generalize the results of \cite{JSGAL} and to associate 
vector bundles to arbitrary equations. This is a quite different mood 
that we will develop in a forthcoming paper (\cite{JSIRR}; meanwhile,
a survey is given in \cite{JSISOMONO}). Here, we give some hints
in remarks \ref{rmk:fibre-functor} and \ref{rmk:vector-bundle}.

Also, let us point out that there has been little effort made
towards systematisation and generalisation. The intent is to
get as efficiently as possible to the striking specific
features of $q$-difference theory. For instance, most of
the results about morphisms between $q$-difference modules
can be obtained by seing these morphisms as meromorphic solutions
of other modules (internal Hom) and  they can therefore be seen
as resulting from more general statements. These facts, evenso
quite often sorites,  deserve to be written. In the
same way, the many regularity properties of the homological
equation $X(qz)A(z)-B(z)X(z)=Y(z)$ should retain some particular
attention and  be clarified in the language of functional
analysis. They are implicitly or explicitly present in many
places in the work of C. Zhang. Last, the $q$-Gevrey filtration
should be translated in terms of factorisation of Stokes
operators, like in \cite{MLR}. \\
\bigskip

Let us now describe the organization of the paper.\\

Notations and conventions  are given in subsection 1.2. \\

Section 2 deals with the recent developments of the theory of
$q$-difference equations and some improvements. In
subsection 2.1, we recall the local classification of fuchsian
systems by means of flat vector bundles as it can be  found in
\cite{JSGAL} and its easy extension to the so-called ``tamely irregular'' 
$q$-difference modules. We then describe the filtration
by the slopes (\cite{JSFIL}). In subsections 2.2, we summarize results 
from \cite{RSZ} about the local
analytic classification of irregular $q$-difference systems,
based on the Stokes sheaf. The lemma \ref{lemma:flatness}
provides a needed improvement about Gevrey decay;
proposition \ref{prop:Rwg=Ka}  and corollary \ref{cor:Rwg=Ka} an
improvement about polynomial normal forms. \\

In chapter 3, we first build our main tool, the algebraic
summation process (theorem \ref{thm:resommation-algebrique}). Its
application to the local classification is then developed in
subsection 3.2. We state there and partially prove the second main
result of this paper (theorem \ref{thm:q-Malgrange-Sibuya}): a
$q$-analog of the Malgrange-Sibuya theorem for the
local analytic classification of linear differential equations.\\

Section 4 is devoted to studying the $q$-Gevrey filtration of the
Stokes sheaf and proving the theorem \ref{thm:q-Malgrange-Sibuya}. 
In subsection 4.1, we show how conditions of flatness (otherwise said, 
of $q$-Gevrey decay) of solutions near
$0$ translate algebraically and how to provide the devissage for
the Stokes sheaf of a ``tamely irregular'' module. In subsection
4.2 we draw some cohomological consequences and we finish the
proof of the theorem \ref{thm:q-Malgrange-Sibuya}. Finally, in
subsection 4.3, we sketch the Stokes sheaf of a general module.\\

The symbol $\Box$  indicates the end of a proof or the absence of
proof if considered straightforward. Theorems,
propositions and lemmas considered as ``prerequisites'' and coming
from the quoted references are not followed by the symbol $\Box$.
\\

\subsection*{Acknowledgements}

The present work 
\footnote{This paper has been submitted (and accepted) 
for publication in the proceedings of the International
Conference in Honor of Jean-Pierre Ramis, held in Toulouse,
september 22-26 2004.}
is directly related with the paper \cite{RSZ},
written in collaboration with Jean-Pierre Ramis and with Changgui
Zhang. It has been a great pleasure to talk with them, confronting
very different points of view and sharing a common excitement. \\

The epigraph at the beginning of this paper is intended to convey
the happiness of wandering and daydreaming in Jean-Pierre Ramis'
garden; and the overwhelming surprise of all its bifurcations.
Like in Borges' story, pathes fork and then unite, the same
landscapes are viewed from many points with renewed pleasure.
This strong feeling of the unity of mathematics without any
uniformity is typical of Jean-Pierre.

% 1.2

\subsection{Notations and general conventions}

We fix once for all a complex number $q \in \C$ such that $|q| > 1$. 
We then define the automorphism $\sq$ on various rings, fields or
spaces of functions by putting $\sq f(z) = f(qz)$. This holds 
in particular for the field $\k$ of complex rational functions,
the ring $\Ra$ of convergent power series and its field 
of fractions $\Ka$, the ring $\Rf$ of formal power series 
and its field of fractions $\Kf$, the ring $\Rwg$ of holomorphic 
germs and the field $\Kwg$ of meromorphic germs in the punctured 
neighborhood of $0$, the ring $\Rw$ of holomorphic functions
and the field $\Kw$ of meromorphic functions on $\C^{*}$; 
this also holds for all modules or spaces of vectors
or matrices over these rings and fields. \\

For any such ring (resp. field) $R$, the $\sq$-invariants elements
make up the subring (resp. subfield) $R^{\sq}$ of constants. For
instance, the field of constants of $\M(\C^{*},0)$ or that of
$\M(\C^{*})$ can be identified with a field of elliptic functions,
the field $\M(\Eq)$ of meromorphic functions over the complex torus
(or elliptic curve) $\Eq = \C^{*}/q^{\Z}$.
We shall use heavily the theta function of Jacobi defined by the
following equality:
$$
\th_{q}(z) = \sum_{n \in \Z} q^{-n(n+1)/2} z^{n}.
$$
This function is holomorphic in $\C^{*}$ with simple zeroes, all
located on the discrete $q$-spiral $[-1;q]$, where we write
$[a;q] = a q^{\Z} \;,\; (a \in \C^{*})$. It satisfies the
functional equation: $\sq \th_{q} = z \th_{q}$. We shall also
use its multiplicative translates $\th_{q,c}(z) = \th_{q}(z/c)$
(for $c \in \C^{*}$); the function $\th_{q,c}$ is holomorphic 
in $\C^{*}$ with simple zeroes, all located on the discrete 
$q$-spiral $[-c;q]$ and satisfies the functional equation: 
$\sq \th_{q,c} = \frac{z}{c} \th_{q,c}$. \\

As is customary for congruence classes,
we shall write $\overline{a} = a \pmod{q^{\Z}}$ for the image
of $a \in \C^{*}$ in the elliptic curve $\Eq$. This notation
extends to a subset $A$ of $\C^{*}$, so that $\overline{A}$
does \emph{not} denote its topological closure. Then, for a divisor
$D = \sum n_{i} [\alpha_{i}]$ over $\Eq$ (\emph{i.e.}, the $n_{i} \in \Z$,
the $\alpha_{i} \in \Eq$), we shall write 
$ev_{\Eq}(D) = \sum n_{i} \alpha_{i} \in \Eq$ for its evaluation, 
computed with the group law on $\Eq$. \\

Let $K$ denote any one of the forementioned fields of functions.
Then, we write $\D = K\l<\sigma,\sigma^{-1}\r>$ for the \"{O}re algebra
of non commutative Laurent polynomials characterized by the relation
$\sigma . f = \sq(f) . \sigma$. We now define the category of 
$q$-difference modules in three clearly equivalent ways:
\begin{eqnarray*}
\DM & = & \{(E,\Phi) \;/\; E \text{~a~} K \text{-vector space of finite rank~},
                        \Phi : E \rightarrow E \; 
                               \text{~a~} \sq \text{-linear} \text{~map}\} \\
    & = & \{(K^{n}, \Phi_{A}) \;/\; A \in GL_{n}(K) ,
                                    \Phi_{A}(X) = A^{-1} \sq X \} \\
    & = & \text{~finite length left~} \D \text{-modules}.
    \end{eqnarray*}
This is a $\C$-linear abelian rigid tensor category, hence
a tannakian category. For basic facts and terminology about these,
see \cite{JSGAL}, \cite{vdPS}, \cite{DM}, \cite{DF}. 
Last, we note that all objects in $\DM$ have the form $\D / \D P$.

%%%%%%%%%%%%%%%%%%%%%%%%%%%%%%%%%%%%%%%%%%%%%%%%%%%%%%%%%%%%%%%%%%%%%%%%%%%%%

% 2

\section{Local analytic classification} 
\label{section:classification}

% 2.1

\subsection{Devissage of irregular equation (\cite{JSGAL},\cite{JSFIL})}
\label{subsection:devissage}

\subsubsection*{Fuchsian and tamely irregular modules}

For a $q$-difference module $M$ over any of the fields $\k$, $\Ka$, 
$\Kf$, it is possible to define its Newton polygon at $0$, or, 
equivalently, the slopes of $M$, which we write in descending 
order : $\mu_{1} > \cdots > \mu_{k} \in \Q$, and their
multiplicities $r_{1},\ldots,r_{k} \in \N^{*}$. 
The module $M$ is said to be
pure of slope $\mu_{1}$ if $k = 1$ and fuchsian if it is pure of
slope $0$. The latter condition is equivalent to $M$ having the shape
$M = (K^{n}, \Phi_{A})$ with $A(0) \in GL_{n}(\C)$. There are also 
criteria of growth (or decay) of solutions near $0$, see further below, 
in section \ref{subsection:irregular-theory}, the subsection about 
flatness conditions. \\

Call $\EE$ the category $DiffMod(\k,\sq)$ of rational equations.  
Fuchsian modules at $0$ and $\infty$ over $\k$ make up a tannakian
subcategory $\EE_{f}$ of $\EE$. In order to study them, one ``localizes'' 
these categories by extending the class of morphisms, precisely, by allowing 
morphisms defined over $\Ka$. This gives ``thickened'' categories $\EE^{(0)}$ 
and $\EE_{f}^{(0)}$. A classical lemma says that any fuchsian system 
is locally equivalent to one with constant coefficients. 
This suggests 
the introduction of the full subcategory $\P_{f}^{(0)}$ of $\EE_{f}^{(0)}$ 
made up of ``flat'' objects, that is, the $(\k^{n}, \Phi_{A})$ with 
$A \in GL_{n}(\C)$. Thus, the inclusion of $\P_{f}^{(0)}$ into 
$\EE_{f}^{(0)}$ is actually an isomorphism of tannakian categories. \\

To any $A \in GL_{n}(\C)$ one associates the holomorphic vector bundle 
$F_{A}$ over $\Eq$ obtained by quotienting $\C^{*} \times \C^{n}$ 
by the equivalence relation $\sim_{A}$ generated by the relations 
$(z,X) \sim_{A} (qz,AX)$. This defines a functor from $\P_{f}^{(0)}$ 
to the category $Fib_{p}(\Eq)$ of flat holomorphic vector bundles over 
$\Eq$. This is an equivalence of tannakian categories. 
Note that the classical lemma alluded above equally holds for any 
fuchsian $q$-difference module over $\Ka$ or over $\Kf$, which implies 
that this local classification applies to $\DMa$ and $\DMf$ as well. 
The galoisian aspects of this local correspondance and its global 
counterpart are detailed in \cite{JSGAL}. \\

A pure module of \emph{integral} slope $\mu$ over $K = \Ka$ or $\Kf$
has the shape $(K^{n}, \Phi_{z^{-\mu} A})$ with $A \in GL_{n}(\C)$.
For such a module, the above construction of a vector bundle extends
trivially, yielding the tensor product of a flat bundle by a line
bundle. We shall call pure such a bundle. \\

Direct sums of pure modules play a special role in \cite{JSFIL},
\cite{RSZ} and in the present paper. We shall call them tamely
irregular, in an intended analogy with tamely ramified extensions
in algebraic number theory: for us, they are irregular objects without
wild monodromy, as follows from \cite{RSZ}. The category of 
tamely irregular modules with integral slopes over $\Ka$ can,
for the same reasons as above, be seen either as a subcategory
of $\DMa$ or of $\EE^{(0)}$. We write it $\EE_{mi,1}^{(0)}$
\footnote{The subscript ``mi'' stands for ``moder\'ement 
irr\'egulier'', the subscript $1$ for restricting to slopes
with denominator $1$.}.
It is generated (as a tannakian category) by the fuchsian modules
and by the pure module $(\Ka,z^{-1} \sq)$ of slope $1$. We can thus
associate to any such module a direct sum of pure modules, thereby
defining a functor from $\EE_{mi,1}^{(0)}$ to the category $Fib(\Eq)$ 
of holomorphic vector bundles over $\Eq$. This functor is easily seen
to be compatible with all linear operations (it is  a functor of 
tannakian categories). \\

\subsubsection*{Filtration by the slopes}

The following is proved in \cite{JSFIL}:

\begin{thm}
Let the letter $K$ stand for the field $\Ka$ (convergent case)
or the field $\Kf$ (formal case). In any case, any object $M$ of $\DM$ 
admits a unique filtration $(F^{\geq \mu}(M))_{\mu \in \Q}$ by subobjects
such that each $F^{(\mu)}(M) = \frac{F^{\geq \mu}(M)}{F^{> \mu}(M)}$ is pure
of slope $\mu$. The $F^{(\mu)}$ are endofunctors of $\DM$ and
$gr = \bigoplus F^{(\mu)}$ is a faithful exact $\C$-linear 
$\otimes$-compatible functor and a retraction of the inclusion 
of $\EE_{mi}^{(0)}$ into $\EE^{(0)}$. In the formal case, $gr$ is isomorphic
to the identity functor.
\end{thm}

\emph{\textbf{From now on, we only consider the full subcategory
$\EE_{1}^{(0)}$ of modules with integral slopes}}.
The notation $\EE_{1}^{(0)}$ will be justified \emph{a posteriori}
by the fact that all its objects are locally equivalent to objects
of $\EE$ (existence of a normal polynomial form).
This is an abelian tensor subcategory of $\EE^{(0)}$ and the functor 
$gr$ retracts $\EE_{1}^{(0)}$ to $\EE_{mi,1}^{(0)}$. 
We also introduce notational conventions which will be used all along 
this paper for a module $M$ in $\EE_{1}^{(0)}$ and its associated graded
module $M_{0} = gr(M)$, an object of $\EE_{mi,1}^{(0)}$. \\

The module $M$ may be given the shape 
$M = (\Ka^{n},\Phi_{A})$, with:
\begin{equation} \label{equation:forme-canonique}
A = A_{U} \underset{def}{=}
\begin{pmatrix}
z^{- \mu_{1}} A_{1}  & \ldots & \ldots & \ldots & \ldots \\
\ldots & \ldots & \ldots  & U_{i,j} & \ldots \\
0      & \ldots & \ldots   & \ldots & \ldots \\
\ldots & 0 & \ldots  & \ldots & \ldots \\
0      & \ldots & 0       & \ldots & z^{- \mu_{k}} A_{k}    
\end{pmatrix},
\end{equation}
where $\mu_{1} > \cdots > \mu_{k}$ are integers,
$r_{i} \in \N^{*}$, $A_{i} \in GL_{r_{i}}(\C)$ ($i = 1,\ldots,k$) and
$$
U = (U_{i,j})_{1 \leq i < j \leq k} \in 
\underset{1 \leq i < j \leq k}{\prod} \Ma_{r_{i},r_{j}}(\Ka).
$$
The associated graded module is then a direct sum
$M_{0} = P_{1} \oplus \cdots \oplus P_{k}$, where,
for $1 \leq i < j \leq k$, the module $P_{i}$ is pure 
of rank $r_{i}$ and slope $\mu_{i}$ and can be put into
the form $P_{i} = (\Ka^{r_{i}},\Phi_{z^{-\mu_{i} A_{i}}})$.
Therefore, one has $M_{0} = (\Ka^{n},\Phi_{A_{0}})$, where
the matrix $A_{0}$ is block-diagonal:
\begin{equation} \label{equation:forme-canonique-diagonale}
A_{0} =
\begin{pmatrix}
z^{- \mu_{1}} A_{1}  & \ldots & \ldots & \ldots & \ldots \\
\ldots & \ldots & \ldots  & 0 & \ldots \\
0      & \ldots & \ldots   & \ldots & \ldots \\
\ldots & 0 & \ldots  & \ldots & \ldots \\
0      & \ldots & 0       & \ldots & z^{- \mu_{k}} A_{k}    
\end{pmatrix}.
\end{equation}

\subsubsection*{The set of analytic isoformal classes}

This section comes from \cite{RSZ}. The definitions here should
be compared to those in \cite{Varadarajan}, p. 29 or \cite{BV}. \\

In $\DMf$, the canonical filtration of a module $M$ is split; 
more precisely, the associated graded module $gr(M)$ is the unique 
formal classifyer of $M$. The isoformal analytic classification is
therefore the same as the isograded classification, whence the
following definitions.

\begin{defn}
Let $P_{1},\ldots,P_{k}$ be pure modules with ranks $r_{1},\ldots,r_{k}$ 
and with integral slopes $\mu_{1} > \cdots > \mu_{k}$. The module
$M_{0} = P_{1} \oplus \cdots \oplus P_{k}$ has rank
$n = r_{1} + \cdots + r_{k}$. We shall write $\F(M_{0})$ for the set
of equivalence classes of pairs $(M,g)$ of a module $M$ and an
isomorphismism $g : \gr(M) \rightarrow M_{0}$, where $(M,g)$ is said to be 
equivalent to $(M',g')$ if there exists a morphism $u : M \rightarrow M'$ 
such that $g = g' \circ \gr(u)$ ($u$ is automatically an isomorphism).
\end{defn}

We write $\G$ for the algebraic subgroup of $GL_{n}$ made up of matrices 
of the form:
\begin{equation} \label{equation:automorphisme}
F = 
\begin{pmatrix}
I_{r_{1}} & \ldots & \ldots & \ldots & \ldots \\
\ldots & \ldots & \ldots  & F_{i,j} & \ldots \\
0      & \ldots & \ldots   & \ldots & \ldots \\
\ldots & 0 & \ldots  & \ldots & \ldots \\
0      & \ldots & 0       & \ldots & I_{r_{k}}     
\end{pmatrix}.
\end{equation}
Its Lie algebra $\g$ consists in matrices of the form:
\begin{equation} \label{equation:endomorphisme}
f = 
\begin{pmatrix}
0_{r_{1}} & \ldots & \ldots & \ldots & \ldots \\
\ldots & \ldots & \ldots  & f_{i,j} & \ldots \\
0      & \ldots & \ldots   & \ldots & \ldots \\
\ldots & 0 & \ldots  & \ldots & \ldots \\
0      & \ldots & 0       & \ldots & 0_{r_{k}}     
\end{pmatrix}.
\end{equation}
For $F$ in $\G$, we shall write $F[A] = \l(\sq F\r) A F^{-1}$ for
the result of the gauge transformation $F$ on the matrix $A$. \\

We shall identify $P_{i}$ with $(\Ka^{r_{i}},\Phi_{z^{-\mu_{i}} A_{i}})$, 
where $A_{i} \in GL_{r_{i}}(\C)$. The datum of a pair $(M,g)$ then
amounts to that of a matrix $A$ in the form \ref{equation:forme-canonique}.
Two such matrices $A,A'$ are equivalent iff there exists a matrix
$F \in \G(\Ka)$ such that $F[A] = A'$. \\

Write
$\G^{A_{0}}(\Kf) = \{F \in \G(\Kf) \;/\; F[A_{0}] \in GL_{n}(\Ka) \}$.
The subgroup $\G(\Ka)$ of $\G(\Kf)$ operates at left on the latter
(by translation) and $\G^{A_{0}}(\Kf)$ is stable for that operation. 
The theory in the previous section entails:
$$
\forall (U_{i,j})_{1 \leq i < j \leq k} \in 
\prod_{1 \leq i < j \leq k} \Ma_{r_{i},r_{j}}(\Ka) \;,\;
\exists ! \hat{F} \in \G(\Kf) \;:\;
\hat{F} [A_{0}] = A_{U}.
$$
This $\hat{F}$ will be written $\hat{F}(U)$. Its existence can also
be proved by direct computation, solving by iteration the fixpoint
equation of the $z$-adically contracting operator: 
$F \mapsto \l(A_{U}\r)^{-1} \l(\sq \hat{F}\r) A_{0}$.
It follows that the unique formal gauge transformation of $\G(\Kf)$
taking $A_{U}$ to $A_{V}$ is $\hat{F}(U,V) = \hat{F}(V) \hat{F}(U)^{-1}$.
Besides, $A_{U}$ is equivalent to $A_{V}$ in the above sense
if and only if $\hat{F}(U,V) \in \G(\Ka)$, or, equivalently,
$\hat{F}(V) \in \G(\Ka) \hat{F}(U)$. 
This translates into the following lemma.

\begin{prop} \label{prop:formel-modulo-convergent}
Sending $A_{U}$ to $\hat{F}(U)$ induces a one-to-one correspondance 
between $\F(M_{0})$ and the
left quotient $\G(\Ka) \backslash \G^{A_{0}}(\Kf)$.
\end{prop}

One thus recognizes
in isoformal classification a classical problem of summation of 
divergent power series. In order to illustrate the possible strategies,
we shall end this section by examining a specific example. We shall try, 
as far as possible, to mimic the methods and the terminology of the 
``classical'' theory (Stokes operators for linear differential equations 
and summation in sectors along directions).

\begin{example} \label{example:exemple-classique}
The module $M_{u} = (\Ka^{2},\Phi_{A_{u}})$ corresponding to the
matrix $A_{u} = \begin{pmatrix} 1 & u \\ 0 & z \end{pmatrix}$ 
is formally isomorphic to its associated graded module $M_{0}$.
More precisely, there exists a formal gauge transformation $F$ 
such that $F[A_{0}] = A_{u}$, that is,$F(qz) A_{0}(z) = A_{u}(z) F(z)$. 
If one moreover requires $F$ to be compatible with the graduation, that 
is, to have the form $F = \begin{pmatrix} 1 & f \\ 0 & 1 \end{pmatrix}$,
then there is unicity of the formal series $f$, which must satisfy the
functional equation:
$$
f(z) = - u(z) + z f(qz).
$$
We call $\hat{f}_{u}$ this unique formal solution (it can be computed 
by iterating the above fixpoint equation) and $\hat{F}_{u}$ the 
corresponding formal gauge transformation. One checks that two such 
matrices $A_{u}$ and $A_{v}$ are analytically equivalent if and only 
if the formal power series $\hat{f}_{u-v} = \hat{f}_{u} - \hat{f}_{v}$ 
is convergent. In this case (two slopes), the problem is additive. \\

For $u = 1$, the unique solution is 
$$
\hat{f}_{1} = - \sum_{n \geq 0}  q^{n(n-1)/2} z^{n},
$$
the so-called Tschakaloff series (up to the sign). It is divergent 
and may be seen as a natural $q$-analog of the Euler series. Thus,
$A_{1}$ is not equivalent to $A_{0}$. \\

In general, we apply the formal $q$-Borel-Ramis transform of level $1$, 
defined by
$$
\B_{q,1} \sum_{n} a_{n} z^{n} = \sum_{n} q^{-n(n-1)/2} a_{n} \xi^{n}.
$$
It sends convergent series to series with an infinite radius of
convergence. Our functional equation is transformed into:
$$
(1-\xi) \B_{q,1}f(\xi) = - \B_{q,1}u(\xi).
$$
The existence of a convergent solution $f$ has only one obstruction,
the number $\nu = \B_{q,1} u (1)$. This number can therefore be considered 
as the unique analytic invariant of $M_{u}$ within the formal class 
of $M_{0}$. It can
also be considered as giving a normal form, since $A_{\nu}$ is the
unique matrix in the analytic class of $A_{u}$ such that $\nu \in \C$.
It is a particular case of normal polynomial form (see further below). \\

The functional equation can also be solved by a variant of the method
of ``varying constants''. We look for the solution in the form
$g = \th_{q,\la} f$. For convenience, we also write $v = \th_{q,\la} u$,
which is an element of $\Rwg$. We compare their Laurent series
coefficientwise and get:
$$
\forall n \in \Z \;,\; (1 - \la q^{n}) g_{n} = v_{n}.
$$
If $\la \not\in [1;q]$ (prohibited direction of summation),
there is a unique solution $g \in \Rwg$ (it does converge
where it should), thus a unique solution $f \in \Kwg$ such that 
$\th_{q,\la} f$ has no poles in $\C^{*}$. We then get a unique
solution $f_{\la,u}$ with (at most) simple poles over $[-\la;q]$:
it is the summation of $\hat{f}_{u}$ in the direction 
$\overline{\la} \in \Eq$ and its ``sector'' of validity is
(the germ at $0$) of $\C^{*} \setminus [-\la;q]$,
the preimage by the canonical projection $\C^{*} \rightarrow \Eq$ 
of the Zariski open set $\Eq \setminus \{\overline{-\la}\}$. \\

There is another way of looking at this summation process, 
with a deeper analytical meaning. We can consider $\B_{q,1} f(\xi)$
as a meromorphic function $\phi$ over the $\xi$-plane and apply to it
some $q$-analog of the Laplace transform. In our case, putting
$$
\mathcal{L}_{q,1}^{\la} \phi(z) =
\sum_{\xi \in [\la;q]} \frac{\phi(\xi)}{\th_{q}(z/\xi)},
$$
gives again $f_{\la,u}$. This \emph{discrete summation} process 
is due to Changgui Zhang see (\cite{Zha02}, also see \cite{DRSZ}) 
and it is heavily used in \cite{RSZ}. In this work, we rather use 
the first more algebraic and more naive method.
\end{example}

% 2.2

\subsection{Classification through the Stokes sheaf 
(\cite{RSZ},\cite{RSZcras1})}
\label{subsection:irregular-theory}

\subsubsection*{The Stokes sheaf and its Lie algebra}

First, we recall the relevant definitions about asymptotic expansions.
The semigroup $\Sigma = q^{-\N}$ operates on $\C^{*}$ with quotient $\Eq$
(its \emph{horizon}); in the classical setting, one would rather have 
an operation of the semigroup $\Sigma = e^{]-\infty,0]}$ with horizon
the circle $S^{1}$ of directions. We consider as sectors the germs
at $0$ of invariant open subsets of $\C^{*}$. We introduce two sheaves
of differential algebras over $\C^{*}$ by putting, for any sector $U$:
\begin{eqnarray*}
\B(U) & = & \{f \in \O(U) \;/\; f 
\text{~is bounded on every invariant relatively compact subset of~}
U\} \\
\A'(U) & = & \{ f \in \O(U) \;/\; \exists \hat{f} \in \Rf \;:\;
\forall n \in \N \;,\; z^{-n}\l(f - S_{n-1} \hat{f}\r) \in \B(U)\},
\end{eqnarray*}
where, as usual,$S_{n-1} \hat{f}$ stands for the truncation.
For any sector $U$, we write $U_{\infty} = U /\Sigma$ for its horizon 
(an open subset of $E_{q}$). We now define a sheaf of differential
algebras over $E_{q}$ by putting:
$$
\A(V) = \underset{\rightarrow}{\lim} \; \A'(U),
$$
the direct limit being taken for the system of those open subsets $U$
such that their horizon is $U_{\infty} = V$. There is a natural
morphism from $\A$ to the constant sheaf with fibre $\Ra$ over $\Eq$
and it is an epimorphism ($q$-analog of Borel-Ritt lemma). We call $\A_{0}$ 
its kernel, the sheaf of infinitely flat functions. For instance,
it is easy to see that a solution of fuchsian equation divided
by a product of theta functions is flat within its domain
(more on this in the next subsection). \\

We then write $\La_{I} = I_{n} + Mat_{n}(\A_{0})$ for the subsheaf
of groups of $GL_{n}(\A)$ made up of matrices infinitely tangent
to the identity and we put $\La_{I}^{\G} = \La_{I} \cap \G(\A)$.
This is a sheaf of matrices of the form \ref{equation:automorphisme}
with all the $F_{i,j}$ flat. 
Last, for a module $M = (\Ka^{n},\Phi_{A})$, we consider the subsheaf
$\La_{I}(M)$ of $\La_{I}^{\G}$ whose sections $F$ satisfy the equality: 
$F[A] = A$ (automorphisms of $M$ infinitely tangent to identity).
This is the \emph{Stokes sheaf} of the module $M$. \\

Note for further use that $\La_{I}(M)$ is a sheaf of unipotent
groups so that one can define algebraically the sheaf $\la_{I}(M)$
of their Lie algebras: we put $\la_{I} = Mat_{n}(\A_{0})$,
$\la_{I}^{\g} = \la_{I} \cap \g(\A)$ (see \ref{equation:endomorphisme})
and take as sections of $\la_{I}(M)$ those sections of $\la_{I}^{\g}$
such that $\l(\sq f\r) A = A f$. Obviously, $f$ is a section of
$\la_{I}(M)$ if and only if $I_{n} + f$ is a section of $\La_{I}(M)$,
or, equivalently, $exp(f)$ is a section of $\La_{I}(M)$. Indeed,
the triangular form and the functional equations are easily
checked, and the flatness properties stem from the well known
fact that, for nilpotent matrices, $f$ and $exp(f) - I_{n}$
are polynomials in each other, without constant terms. 

\subsubsection*{The $q$-analogs of Malgrange-Sibuya theorems}

One showed in \cite{RSZ} the following $q$-analogs of classical 
theorems by Malgrange-Sibuya:

\begin{thm}
There are natural bijective mappings:
$$
\G(\Ra) \diagdown \G^{A_{0}}(\Rf) \rightarrow
\G(\Ka) \diagdown \G^{A_{0}}(\Kf) \rightarrow
H^{1}(E_{q},\La_{I}^{\G}).
$$
\end{thm}

Actually, the following more general theoren is proven in 
\emph{loc. cit.}, dealing with an arbitrary algebraic subgroup 
$G$ of $GL_{n}$. Its proof relies on some heavy analysis 
(Newlander-Nirenberg theorem). 

\begin{thm}
Let $M_{0}$ be as above. There are natural bijective mappings:
$$
\F(M_{0}) \rightarrow
G(\Ra) \diagdown G^{A_{0}}(\Rf) \rightarrow
G(\Ka) \diagdown G^{A_{0}}(\Kf) \rightarrow
H^{1}(E_{q},\La_{I}^{G}),
$$
where $\La_{I}^{G} = \La_{I} \cap G(\A)$.
\end{thm}

The former theorem is deduced from the latter together with the existence 
of asymptotic solutions. One can explicitly build, by discrete 
resummation, privileged cocycles associated to a class in $\F(M_{0})$ 
and to ``Stokes directions''. In the next chapter, I shall exhibit
an algebraic variant of this construction. Morally, it is possible
because the sheaf $\La_{I}(M_{0})$ is almost a vector bundle over
the elliptic curve $\Eq$.

\subsubsection*{Flatness conditions}

Details about the contents of this section can be found 
in \cite{RZ} and \cite{RSZ}; see also the older references 
\cite{RamisGrowth} and \cite{RamisJPRTraum}. \\

The above notion of flatness can be refined, introducing 
$q$-Gevrey levels. These may be characterized either in
terms of growth (or decay) of functions near $0$, or in terms of growth
of coefficients of power series. We shall use here the following
simple terminology and facts. \\

We start from a proper germ of $q^{-\N}$ invariant subset $U$
of $(\C^{*},0)$. Then any solution of a fuchsian system that is
holomorphic on $U$ has polynomial growth at $0$ (see for instance 
\cite{JSGAL}); this is for instance true for a quotient of theta 
functions. We say that $f \in \O(U)$ has level of flatness $\geq t$ 
(where $t$ is an integer)
if, for one (hence any) theta function $\th = \th_{q,\la}$, the
function $f |\th|^{t}$ has polynomial growth near $0$. 
We easily get the following implications.

\begin{lemma} \label{lemma:flatness}
(i) For $t > 0$, $t$-flatness implies flatness in the sense
of asymptotics. \\
(ii) Solutions of pure systems of slope $\mu$ are $\mu$-flat. \\
(iii) If a solution of a pure system of slope $\mu$ is $t$-flat
with $t > 0$, then it is $0$.
\end{lemma}

\hfill $\Box$

\subsubsection*{Normal polynomial forms}

The computations will follow the same pattern as in \cite{RSZ},
\cite{RSZcras1}. However, we shall need a slightly more general 
version afterwards (proposition \ref{prop:Rwg=Ka}). \\

We start with a computation with two slopes. Take integers 
$\mu > \mu'$, square invertible matrices
$A \in GL_{r}(\C)$ and $A' \in GL_{r'}(\C)$. 
Just for this section, call $\VV(r,r',\mu,\mu')$
the subspace of $\M_{r,r'}(\Ka)$ spanned by matrices
all of whose coefficients belong to 
$\underset{\mu' \leq k < \mu}{\sum} \C z^{k}$. \\

For $U \in \M_{r,r'}(\Ka)$, write 
$B_{U} = \begin{pmatrix} z^{-\mu} A & U \\ 0 & z^{-\mu'} A' \end{pmatrix}$.
Then, for any such $U$, there exists a unique pair $(F,V)$ with 
$F \in \M_{r,r'}(\Ka)$ and $V \in \VV(r,r',\mu,\mu')$ 
such that the matrix
$\begin{pmatrix} I_{r} & F \\ 0 & I_{r'} \end{pmatrix}$
defines an isomorphism from $B_{U}$ to $B_{V}$.
This amounts to solving:
\begin{equation} \label{equation:equation-homologique}
\l(\sq F \r) (z^{-\mu'} A') - (z^{-\mu} A) F = V - U.
\end{equation}
Successive reductions boil the problem down to example
\ref{example:exemple-classique}. 
We shall write $Red(\mu,A,\mu',A',U)$ for the pair $(F,V)$. \\

Now, we come back to our usual notations \ref{equation:forme-canonique}
and \ref{equation:forme-canonique-diagonale}. We consider the matrix
$A_{U}$ associated to
$U = (U_{i,j}) \in 
\underset{1 \leq i < j \leq k}{\prod} \M_{r_{i},r_{j}}(\Ka)$.
Then, there is a unique pair $(\underline{F},V)$ with
$\underline{F} = (F_{i,j}) \in 
\underset{1 \leq i < j \leq k}{\prod} \M_{r_{i},r_{j}}(\Ka)$
and
$V = (V_{i,j}) \in 
\underset{1 \leq i < j \leq k}{\prod} \VV(r_{i},r_{j},\mu_{i},\mu_{j})$
such that the associated gauge transformation $F \in \G(\Ka)$   
defines an isomorphism from $A_{U}$ to $A_{V}$. The pair $(\underline{F},V)$
can be computed by solving iteratively a system of equations 
of the type \ref{equation:equation-homologique} for $1 \leq i < j \leq k$.
This is done by inductively with the help of the formula:
$$
(F_{i,j},V_{i,j}) = Red\l(\mu_{i},A_{i},\mu_{j},A_{j},
U_{i,j} + \underset{{i < l < j}}{\sum} \l(\sq F_{i,l}\r) U_{l,j}
- \underset{i < l < j}{\sum} V_{i,l} F_{l,j}\r).
$$
What we get is, in essence, the canonical form of Birkhoff and
Guenther. Standing alone, this statement confirms our earlier 
contention in section \ref{subsection:devissage}, to the effect 
that all objects of $\EE_{1}^{(0)}$ are locally equivalent to 
objects of $\EE$. \\

Now, we shall have use for an extension of these results allowing 
for coefficients in $\Rwg$ (instead of $\Ka$).

\begin{prop} \label{prop:Rwg=Ka}
Let $A_{U}$ be as above, but with 
$U = (U_{i,j}) \in 
\underset{1 \leq i < j \leq k}{\prod} \M_{r_{i},r_{j}}(\Rwg)$.
Then, there exists a unique pair $(\underline{F},V)$ with
$\underline{F} = (F_{i,j}) \in 
\underset{1 \leq i < j \leq k}{\prod} \M_{r_{i},r_{j}}(\Rwg)$
and
$V = (V_{i,j}) \in 
\underset{1 \leq i < j \leq k}{\prod} \VV(r_{i},r_{j},\mu_{i},\mu_{j})$
such that the associated gauge transformation $F \in \G(\Rwg)$   
defines an isomorphism from $A_{U}$ to $A_{V}$.
\end{prop}

\Pr 
The same induction as before can be used, and the proof boils down to
the following lemma.

\begin{lemma}
Let $\mu > \mu'$ in $\Z$, $A \in GL_{r}(\C)$,$A' \in GL_{r'}(\C)$
and $U \in \M_{r,r'}(\Rwg)$.
There exists a unique pair $(F,V)$ with 
$F \in \M_{r,r'}(\Rwg)$ and $V \in \VV(r,r',\mu,\mu')$ 
satisfying \ref{equation:equation-homologique}.
\end{lemma}

\Pr 
The same reductions as in \emph{loc. cit.} entail that we may as well 
assume from the beginning that $\mu = 0$ and $\mu' = -1$. 
The equation as written has unknown $F$ and right hand side $V - U$
in a space of rectangular matrices. Call $s$ the rank of this space 
and call $B$ the matrix of its automorphism $F \mapsto A F {A'}^{-1}$ 
relative to some basis. 
Multiplying both sides of \ref{equation:equation-homologique}
by ${A'}^{-1}$ at right, we get an equivalent equation of the shape
$z \; \sq X - B X = Y - Y^{(0)}$, for which we want to show that, 
for arbitrary $B \in GL_{s}(\C)$ and $Y \in \Rwg^{s}$, there is 
a unique solution $(X,Y^{(0)}) \in \Rwg^{s} \times \C^{s}$. 
Note that, replacing
$X = \sum_{n \in \Z} X_{n} z^{n}$, $Y = \sum_{n \in \Z} Y_{n} z^{n}$
and $Y^{(0)}$ respectively by $\sum_{n \in \Z} B^{n} X_{n} z^{n}$, 
$\sum_{n \in \Z} B^{n-1} Y_{n} z^{n}$ and $B^{-1}Y^{(0)}$,
we do not change the conditions on $X$, $Y$, $Y^{(0)}$, and we
are led to study a similar problem with $B = I_{s}$. The latter
problem can be tackled componentwise: we are to show that, for
any $u \in \Rwg^{s}$, there is a unique pair 
$(f,\nu) \in  \Rwg \times \C$ such that $z \sq f - f = u - \nu$
(compare to example \ref{example:exemple-classique}). \\

We apply the $q$-Borel-Ramis transform of level $1$. This clearly sends
$\Rwg$ to $\Rw$: indeed, for any $A > 0$, $A^{n} q^{-n(n-1)/2}$ tends
to $0$ when $n \to \pm \infty$. From the computations in example
\ref{example:exemple-classique}, we deduce that we have to take
$\nu = \B_{q,1} u (1) = \sum_{n \in \Z} q^{-n(n-1)/2} u_{n}$;
we must then prove the existence and unicity of $f$. Replacing $u$
by $u - \nu$, we may assume that $\B_{q,1} u (1) = 0$. 
We write $f'_{n} = q^{-n(n-1)/2} f_{n}$ and $u'_{n} = q^{-n(n-1)/2} u_{n}$
the coefficients of the $q$-Borel-Ramis transforms $\B_{q,1} f$ and
$\B_{q,1} u$. We know that $\sum_{-\infty}^{+\infty} u'_{k} = 0$ and
require that $\forall n \in \Z, f'_{n-1} - f'_{n} = u'_{n}$.
The only possibility allowing $f'_{n} \to 0$ for $n \to \pm \infty$
is given by the two equivalent definitions:
\begin{eqnarray*}
f'_{n} & \underset{def}{=} & \sum_{k = n+1}^{+ \infty} u'_{k} \\
       & \underset{def}{=} & - \sum_{k = -\infty}^{n} u'_{k}.
\end{eqnarray*}
For $n \to + \infty$, we thus take (using the first definition of $f'_{n}$):
$$
f_{n} = q^{n(n-1)/2} \sum_{k = n+1}^{+ \infty} \frac{u_{k}}{q^{k(k-1)/2}}.
$$
By assumption on $u$, there exists $A > 0$ and $C > 0$ such that,
$\forall n \geq 0$, $\left|u_{n}\right| \leq C A^{n}$. Then:
$$
\left|f_{n}\right| \leq \frac{C A^{n+1}}{\left|q\right|^{n}}
\left(
1 + \frac{A}{\left|q\right|^{n+1}} + 
\frac{A^{2}}{\left|q\right|^{(n+1) + (n+2)}} + \cdots
\right),
$$
whence 
$\left|f_{n}\right| = 
\text{O}\left(
\left(\frac{A}{\left|q\right|}\right)^{n}
\right)$
when $n \to + \infty$. \\

On the side of negative powers, putting, for convenience,
$g_{n} = f_{-n}$ and $v_{k} = u_{-k}$ and, using the second definition
for $f'_{n}$, we see that 
$$
g_{n} = q^{n(n+1)/2} \sum_{k = n}^{+ \infty} \frac{v_{k}}{q^{k(k+1)/2}}.
$$
By assumption on $u$, we have, for \emph{any} $B > 0$, 
$\left|v_{k}\right| = \text{O}(B^{k})$ when $k \to + \infty$ and
a similar computation as before then yields that, for \emph{any} $B > 0$, 
$\left|g_{n}\right| = \text{O}(B^{n})$ when $n \to + \infty$,
allowing one to conclude that $f \in \Rwg$ as desired.
\hfill $\Box$ \\

We shall actually need only the following consequence of the proposition.

\begin{cor} \label{cor:Rwg=Ka}
Let $A = A_{U}$ in the canonical form \ref{equation:forme-canonique},
with 
$U = (U_{i,j}) \in 
\underset{1 \leq i < j \leq k}{\prod} \M_{r_{i},r_{j}}(\Rwg)$.
Then, there exists $F \in \G(\Rwg)$ such that $A_{V} = F[A_{U}]$
has the same form, but with 
$V = (V_{i,j}) \in 
\underset{1 \leq i < j \leq k}{\prod} \M_{r_{i},r_{j}}(\Ka)$.
\end{cor}
\hfill $\Box$ \\

Obviously the same properties hold if one replaces $\Rwg$ by $\Rw$.

%%%%%%%%%%%%%%%%%%%%%%%%%%%%%%%%%%%%%%%%%%%%%%%%%%%%%%%%%%%%%%%%%%%%%%%%%%%%%

% 3

\section{Algebraic summation}

% 3.1

\subsection{The algorithm}

We keep the notations $M_{0}$, $A_{0}$, $A_{U}$ of 
\ref{equation:forme-canonique} and 
\ref{equation:forme-canonique-diagonale}
and the corresponding conventions from section \ref{subsection:devissage}.
Also, we shall, for $1 \leq i < j \leq k$, use the abreviation:
$\mu_{i,j} = \mu_{i} - \mu_{j} \in \N^{*}$.

\begin{defn}
(i) A summation divisor adapted to $A_{0}$ 
is a family $(D_{i,j})_{1 \leq i < j \leq k}$ 
of effective divisors over the elliptic curve $\Eq$, each $D_{i,j}$
having degree $\mu_{i,j}$, the family satisfying moreover 
the following compatibility condition:
$$
\forall i,l,j \text{~such that~} 1 \leq i < l < j \leq k \;,\;
D_{i,j} = D_{i,l} + D_{l,j}.
$$
Obviously, it amounts to the same thing to give only the $k-1$
divisors $D_{i,i+1}$, $i = 1,\ldots,k-1$. \\
(ii) We say that the adapted summation divisor 
$(D_{i,j})_{1 \leq i < j \leq k}$ 
is allowed if it satisfies the following conditions:
$$
\forall i,j \text{~such that~} 1 \leq i < j \leq k \;,\;
ev_{\Eq}(D_{i,j}) \not\in 
\overline{(-1)^{\mu_{i,j}} \frac{Sp(A_{i})}{Sp(A_{j})}}.
$$
Here, for $S,T \subset \C^{*}$, we put
$\frac{S}{T} = \{\frac{s}{t} \;/\; s \in S , t \in T\}$;
$\overline{X}$ and $ev_{\Eq}$ were defined in the introduction.
\end{defn}
Note that, for an adapted summation divisor, the condition
of being allowed is a generic one. 

\begin{example}
A special case is that of an adapted summation divisor concentrated 
on a point $\alpha \in \Eq$, that is, each $D_{i,j} = \mu_{i,j} [\alpha]$. 
Then the condition that $D$ is allowed is equivalent to:
$$
\forall i,j \text{~such that~} 1 \leq i < j \leq k \;,\;
\mu_{i,j} \; \alpha \not\in \overline{Sp(A_{i})} - \overline{Sp(A_{j})}.
$$
It is generically (that is, over a non empty Zariski open subset)
satisfied by $\alpha \in \Eq$. 
\end{example}

Now, let $(D_{i,j})_{1 \leq i < j \leq k}$ be a summation divisor
adapted to $A_{0}$. We choose points $a_{l} \in \C^{*}$ for 
$\mu_{k} < l \leq \mu_{1}$ such that, for $1 \leq i < j \leq k$, 
$$
D_{i,j} = \sum_{\mu_{j} < l \leq \mu_{i}} [\overline{a_{l}}].
$$
These certainly exist. We then put:
$$
t_{i} = \th_{q}^{\mu_{k}} \prod_{\mu_{l} < l \leq \mu_{i}} \th_{q,-a_{l}}.
$$

\begin{lemma}
(i) The functions $t_{1},\ldots,t_{k} \in \Kw$ are such that: \\
(i1) For $i = 1,\ldots,k$, $\sq t_{i} = \alpha_{i} z^{\mu_{i}} t_{i}$,
where $\alpha_{i} \in \C^{*}$. \\
(i2) For $1 \leq i < j \leq k$, 
$div_{\Eq}(t_{i}) - div_{\Eq}(t_{j}) = D_{i,j}$
(the notation is explained in the course of the proof). \\
(i3) For $1 \leq i < j \leq k$,
the function $t_{i,j} = \frac{\sq t_{i}}{t_{j}}$ belongs to
$\mathcal{O}(\C^{*})$. \\
(ii)If the summation divisor $(D_{i,j})_{1 \leq i < j \leq k}$ 
is moreover allowed, for $1 \leq i < j \leq k$, the spectra of 
$\alpha_{i} A_{i}$ and $\alpha_{j} A_{j}$ have empty intersection 
on $\Eq$:
$$
\overline{Sp(\alpha_{i} A_{i})} \cap \overline{Sp(\alpha_{j} A_{j})} =
\emptyset.
$$
\end{lemma}

\Pr 

It is an immediate consequence of the properties recalled in the
introduction that these functions $t_{i}$ indeed satisfy (i1).
Moreover, the functional equation implies that the divisor
$div_{\C^{*}}(t_{i})$ of zeroes and poles of $t_{i}$ on $\C^{*}$
is invariant under the action of $q^{\Z}$, so that it makes sense
to consider it as a divisor $\div_{\Eq}(t_{i})$ on $\Eq$
(alternatively, one can consider $t_{i}$ as a section of a line
bundle over $\Eq$ and the notation is then classical). Again
because of the properties of theta-functions, one clearly gets (i2).
Assertion (i3) comes from the equalities:
\begin{eqnarray*}
t_{i,j} & = & \frac{\sq t_{i}}{t_{i}} \times \frac{t_{i}}{t_{j}} \\ 
        & = & \alpha_{i} z^{\mu_{i}} \times 
\text{~a function with positive divisor~}.
\end{eqnarray*}
The function $\th_{q}^{\mu_{i}} \frac{\th_{q}}{\th{q,\alpha_{i}}}$ 
satisfies the same functional equation as $t_{i}$, which means that
their quotient is elliptic so that its divisor on $\Eq$ has trivial
evaluation. Therefore:
$$
div_{\Eq}(t_{i}) = \overline{\frac{(-1)^{\mu_{i}}}{\alpha_{i}}}.
$$
The conclusion (ii) then follows from the definition of an allowed divisor.
\hfill $\Box$ \\

We now introduce a temporary and slightly ambiguous notation.
For an adapted summation divisor $D = (D_{i,j})_{1 \leq i < j \leq k}$,
we write $\Th_{D}$ for the following block-diagonal matrix:
$$
\Th_{D} = 
\begin{pmatrix}
t_{1} I_{r_{1}}  & \ldots & \ldots & \ldots & \ldots \\
\ldots & \ldots & \ldots  & 0 & \ldots \\
0      & \ldots & \ldots   & \ldots & \ldots \\
\ldots & 0 & \ldots  & \ldots & \ldots \\
0      & \ldots & 0       & \ldots & t_{k} I_{r_{k}}    
\end{pmatrix}.
$$
Of course, it does not only depend on $D$, but on a particular
choice of the functions $t_{1},\ldots,t_{k}$ whose existence
has just been established. However, the summation process
that we are defining will produce a result that only depends on $D$.
For a family of recangular blocks $U'_{I,j}$, we shall use the 
following abreviation:
$$
A'_{U'} =
\begin{pmatrix}
A'_{1}  & \ldots & \ldots & \ldots & \ldots \\
\ldots & \ldots & \ldots  & U'_{i,j} & \ldots \\
0      & \ldots & \ldots   & \ldots & \ldots \\
\ldots & 0 & \ldots  & \ldots & \ldots \\
0      & \ldots & 0       & \ldots & A'_{k}    
\end{pmatrix},
$$

\begin{lemma}
(i) The effect of the gauge transformation $\Th_{D}$ is to 
"regularize" the diagonal blocks of $A = A_{U}$:
$\Th_{D}[A] = A'_{U'}$, where, for $1 \leq i < j \leq k$, 
$U'_{i,j} = t_{i,j} U_{i,j} \in M_{r_{i},r_{j}}(\Rwg)$
and, for $i = 1,\ldots,k$, $A'_{i} = \alpha_{i} A_{i} \in GL_{r_{i}}(\C)$.
If moreover the adapted summation divisor $D$ is allowed, then,
for $1 \leq i < j \leq k$, 
$\overline{Sp(A'_{i})} \cap \overline{Sp(A'_{j})} = \emptyset$. \\
(ii) Suppose we started with $A_{U}$ in polynomial normal form.
Then we get $A'_{U'}$ such that $U'_{i,j} \in M_{r_{i},r_{j}}(\Rw)$.
\end{lemma}

\Pr 
The computations are immediate.
\hfill $\Box$ \\

We shall now take two matrices $A_{U}$ and $A_{V}$ in the formal
class of $A_{0}$, flatten their slopes through the gauge transformation
$\Th_{D}$, and then link the resulting matrices $A'_{U'}$ and $A'_{V'}$ 
by an isomorphism defined over $\C^{*}$. This relies on the following

\begin{prop}
(i) Let
$$
A'_{U'} = \begin{pmatrix}
A'_{1}  & \ldots & \ldots & \ldots & \ldots \\
\ldots & \ldots & \ldots  & U'_{i,j} & \ldots \\
0      & \ldots & \ldots   & \ldots & \ldots \\
\ldots & 0 & \ldots  & \ldots & \ldots \\
0      & \ldots & 0       & \ldots & A'_{k}    
\end{pmatrix}
\quad \text{~and~} \quad
A'_{V'} = \begin{pmatrix}
A'_{1}  & \ldots & \ldots & \ldots & \ldots \\
\ldots & \ldots & \ldots  & V'_{i,j} & \ldots \\
0      & \ldots & \ldots   & \ldots & \ldots \\
\ldots & 0 & \ldots  & \ldots & \ldots \\
0      & \ldots & 0       & \ldots & A'_{k}    
\end{pmatrix},
$$
where, for $i = 1,\ldots,k$, $A'_{i} \in GL_{r_{i}}(\C)$ are such that,
for $1 \leq i < j \leq k$, 
$\overline{Sp(A'_{i})} \cap \overline{Sp(A'_{j})} = \emptyset$
and, for $1 \leq i < j \leq k$, 
$U'_{i,j}, V'_{i,j} \in M_{r_{i},r_{j}}(\Rwg)$.
Then, there exists a unique $F' \in \G(\Rwg)$ such that
$F'[A'_{U'}] = A'_{V'}$. \\
(ii) If , for $1 \leq i < j \leq k$, 
$U'_{i,j}, V'_{i,j} \in M_{r_{i},r_{j}}(\Rw)$,
then $F' \in \G(\Rw)$
\end{prop}

\Pr 
We have to solve inductively the system of equations:
$$
\l(\sq F'_{i,j}\r) - A'_{i} F'_{i,j} =
V'_{i,j} - U'_{i,j} +
\sum_{i < l < j} V'_{i,l} F'_{l,j} - 
\sum_{i < l < j} \l(\sq F'_{i,l}\r) U'_{l,j}.
$$
The induction is the same as the one we met when buliding
normal polynomial forms. The proposition then 
follows from the following lemma.

\begin{lemma}
(i) Let $B \in GL_{s}(\C)$ and $C \in GL_{t}(\C)$ be invertible
complex matrices such that 
$\overline{Sp(B)} \cap \overline{Sp(C)} = \emptyset$.
Then, for $Y' \in M_{t,s}(\Rwg)$, the equation:
$$
\l(\sq X'\r) B - C X' = Y'
$$
has a unique solution $X' \in M_{t,s}(\Rwg)$. \\
(ii) If $Y' \in M_{t,s}(\Rw)$, then $X' \in M_{t,s}(\Rw)$.
\end{lemma}

\Pr 
We write the Laurent series:
$$
X' = \sum_{n \in \Z} X'_{n} z^{n} \quad \text{and} \quad
Y' = \sum_{n \in \Z} Y'_{n} z^{n}.
$$
By identification, we obtain $X'_{n} = \Phi_{q^{n} B,C}^{-1}(Y'_{n})$,
where $\Phi_{q^{n} B,C}$ is the automorphism $M \mapsto M (q^{n} B) - C M$ 
of $M_{t,s}(\C)$; that it is indeed an automorphism comes from the
assumption that $q^{n} B$ and $C$ have non intersecting spectra.
For $n \to + \infty$, 
$\Phi_{q^{n} B,C}^{-1} \sim q^{-n} \Phi_{B,0}^{-1}$
and, for 
$n \to - \infty$, $\Phi_{q^{n} B,C}^{-1} \to \Phi_{0,C}^{-1}$.
Taking $]r,R[ \times e^{\imath \R}$ to be the annulus of convergence
of $Y'$, we conclude that the annulus of convergence of $X'$ is 
$]r,|q| R[ \times e^{\imath \R}$.
Annuli of definition actually grow, again an illustration of the
good regularity properties of the homological equation.
\hfill $\Box$ \\

Putting it all together, we now get our first fundamental theorem.

\begin{thm} \label{thm:resommation-algebrique}
(i) Let $A_{U},A_{V}$ be defined as above, in the formal class of $A_{0}$.
Then, there exists a unique $F \in \G(\Kwg)$ such that $F[A_{U}] = A_{V}$ 
and, for $1 \leq i < j \leq k$, $\div_{\Eq}(F_{i,j}) \geq - D_{i,j}$
(the notation is explained in the course of the proof). \\
(ii) If $A_{U},A_{V}$ are in polynomial normal form, $F \in \G(\Kw)$.
\end{thm}

\Pr 
We put $A'_{U'} = \Th_{D}[A_{U}]$ and $A'_{V'} = \Th_{D}[A_{V}]$,
then $F[A_{U}] = A_{V}$ is equivalent to $F'[A'_{U'}] = A'_{V'}$,
where $F' = \Th_{D} F \Th_{D}^{-1}$. The matrices $F$ and $F'$
together are upper-triangular with diagonal blocks 
$I_{r_{1}},\ldots,I_{r_{k}}$ and their over-diagonal blocks
are related by the relations: $F_{i,j} = \frac{t_{j}}{t_{i}}F'_{i,j}$.
This implies the unicity of $F \in \G(\Kwg)$ (resp. $\G(\Kw)$)
subject to the constraint that the coefficients of $F_{i,j}$ belong 
to $\frac{t_{j}}{t_{i}} \Rwg$ (resp. $\frac{t_{j}}{t_{i}} \Rw$). 
Since $div_{\Eq}(t_{i}) - div_{\Eq}(t_{j}) = D_{i,j}$, this proves
(and explains) the given condition. 
\hfill $\Box$ \\

As a matter of notation, we shall write $F_{D}(U,V)$ for the $F$
obtained in the theorem: it does indeed depend solely on $D$. We
see it as the canonical resummation of $\hat{F}(U,V)$ along the
``direction'' $D$. We shall write in particular $F_{D}(U) = F_{D}(0,U)$. \\

Let us call $\G_{D}(\Kwg)$ (resp. $\G_{D}(\Kw)$) the subset of
$\G(\Kwg)$ (resp. of $\G(\Kw)$) defined by the constraints:
for $1 \leq i < j \leq k$, $\div_{\Eq}(F_{i,j}) \geq - D_{i,j}$.

\begin{cor}
$F_{D}(U,V) = F_{D}(V) F_{D}(U)^{-1}$.
\end{cor}

\Pr
Actually, $\G_{D}(\Kwg) = \Th_{D} \G(\Rwg) \Th_{D}^{-1}$ and
$\G_{D}(\Kw) = \Th_{D} \G(\Rw) \Th_{D}^{-1}$, so that these subsets
are subgroups. Then the statement follows from the unicity
property in the theorem.
\hfill $\Box$

\begin{cor}
The conclusion of the theorem still holds if one only assumes that
$U = (U_{i,j}) \in 
\underset{1 \leq i < j \leq k}{\prod} \M_{r_{i},r_{j}}(\Rwg)$.
\end{cor}

\Pr
This immediately follows from corollary \ref{cor:Rwg=Ka}.
\hfill $\Box$

\begin{rmk}
It is not difficult to prove that $\hat{F}(U,V)$ is the aymptotic
expansion of $F_{D}(U,V)$ in the sense of section
\ref{subsection:irregular-theory}. One first has to extend 
the definitions so as to allow for a pole at $0$.The proof then
proceeds in two steps.
\begin{enumerate}
\item{First, one proves that, in its domain of definition,
$F_{D}(U,V)$ is a section of the sheaf $z^{d} \B$ for some
$d \in \Z$. This is done using only the functional equation
that it satisfies, and studying inductively its upper diagonal
blocks $F_{i,j}$.}
\item{Then, one proves that the operator 
$F \mapsto A_{V}^{-1} \l(\sq F\r) A_{U}$, 
sends $\G(z^{d} \B)$ to $\G(z^{d+1} \B)$.
Starting from $F_{D}(U,V)$ and iterating yields the conclusion.}
\end{enumerate}
Actually, in \cite{RSZ}, a stronger result is proved. It relies 
on a refined definition of asymptotics taking in account the position 
of poles; this is essential to get summation by discrete integral
formulas.
\end{rmk}

\begin{rmk} \label{rmk:fibre-functor}
To give our theorem its functorial meaning, one should proceed 
as follows. One generalizes the construction of a vector bundle
$F_{M}$ from a $q$-difference module $M$. This defines a fibre functor
$\omega$ over $\Eq$. Then, for each $D$, if one restricts to an appropriate 
subcategory of $\EE^{(0)}_{1}$, $M \leadsto F_{D}(M)$ is an
isomorphism from the fibre functor $\omega \circ gr$ to $\omega$. 
On the other hand, endowing $F_{M}$ with the filtration coming from 
that of $M$, one defines an enriched functor and the underlying
principle of all our uses of the homological equation is that 
this functor is fully faithful. This is exploited in \cite{JSIRR}.
\end{rmk}

% 3.2

\subsection{Applications to classification}

\subsubsection*{One direction of summation}

Let $A_{U},A_{V}$ be defined over $\Ka$. Suppose $A_{U}$ and $A_{V}$ 
are analytically equivalent. Then the power series $\hat{F}(U,V)$ 
is convergent and satisfies the conclusion of theorem
\ref{thm:resommation-algebrique}, so that, 
by unicity, $F_{D}(U,V) = \hat{F}(U,V)$ for any allowed summation divisor.
Conversely:

\begin{prop} \label{prop:tout-dans-les-poles} 
Suppose $F_{D}(U,V) \in \G(\Rwg)$. Then $A_{U}$ and $A_{V}$ are analytically 
equivalent (and all the above holds).
\end{prop}

\Pr 
The gauge transform $F = F_{D}(U,V)$ is obtained by solving the system
of equations:
$$
z^{-\mu_{j}} \l(\sq F_{i,j} \r) A_{j} - z^{-\mu_{i}} A_{i} F_{i,j} =
V_{i,j} - U_{i,j} + \sum_{i < l < j} V_{i,l} F_{l,j} - 
\sum_{i < l < j} \l(\sq F_{i,l}\r) U_{l,j}.
$$
By induction, we are reduced to the following lemma:

\begin{lemma}
Let $\mu > \mu'$ in $\Z$, $A \in GL_{r}(\C)$,$A' \in GL_{r'}(\C)$
and $Y \in \M_{r,r'}(\Ka)$. Let $X \in \M_{r,r'}(\Rwg)$ be a solution
of the equation:
$$
\l(\sq X \r) (z^{-\mu'} A') - (z^{-\mu} A) X = Y.
$$
Then, one actually has $X \in \M_{r,r'}(\Ka)$.
\end{lemma}

\Pr 
Going to the Laurent series and identifying coefficients,
one finds:
$$
\forall n \in \Z \;,\; 
q^{n+\mu'} X_{n + \mu'} A' - A X_{n+\mu} = Y_{n}.
$$
Since $Y \in \M_{r,r'}(\Ka)$, $Y_{n} = 0$ for $n << 0$.
Therefore, for $n << 0$, writing $d = \mu - \mu' \in \N^{*}$,
one has $X_{n} = q^{-n} A X_{n+d} {A'}^{-1}.$ Since $|q| > 1$,
either $X_{n} = 0$ for $n << 0$, or the coefficients of $X$ are
rapidly growing for indices near $- \infty$ prohibiting convergence
and contradicting the assumption that $X \in \M_{r,r'}(\Rwg)$.
\hfill $\Box$ \\

In order to make this a statement about classification, we introduce
one more notation. We write:
$$
\G_{D}^{A_{0}}(\Kwg) = \{F \in \G_{D}(\Kwg) \;/\; F[A_{0}] \in GL_{n}(\Ka)\}.
$$
Clearly, the subset $\G_{D}^{A_{0}}(\Kwg)$ of the group $\G_{D}(\Kwg)$
is stable under the action by left translations of the subgroup $\G(\Rwg)$.
Now, the above proposition immediately entails:

\begin{prop}
Mapping $A_{U}$ to $F_{D}(U)$ yields a bijection:
$$
\F(M_{0}) \rightarrow \G(\Rw) \diagdown \G_{D}^{A_{0}}(\Kwg).
$$
\end{prop}
\hfill $\Box$ \\

This is strikingly similar to the corresponding ``formal modulo
analytic'' description in proposition \ref{prop:formel-modulo-convergent}.

\subsubsection*{Varying the direction of summation}

Let $D = (D_{i,j})_{1 \leq i < j \leq k}$ be an allowed summation 
divisor for $M_{0}$, $A_{0}$. We consider as its support and write
$Supp(D)$ the union $\underset{1 \leq i < j \leq k}{\bigcup} Supp(D_{i,j})$
and define the following Zariski open subset of $\Eq$:
$V_{D} = \Eq \setminus \text{Supp}(D)$. We also write $U_{D}$ for 
the preimage of $V_{D}$ in $\C^{*}$. Thus, the elements of $\G_{D}(\Kwg)$ are 
holomorphic germs over $(U_{D},0)$. \\

Now let $D'$ be another allowed summation divisor. Then, for any
$A_{U}$ in the formal class of $A_{0}$, the gauge transformation
$F_{D,D'}(U) \underset{def}{=} F_{D}(U)^{-1} F_{D'}(U)$ sends $A_{0}$ 
to itself. It is holomorphic on the open subset $(U_{D} \cap U_{D'},0)$. 
We call $U_{D,D'}$ the sector $U_{D} \cap U_{D'}$, which is the preimage 
of the open subset $V_{D,D'} = V_{D} \cap V_{D'}$ of $\Eq$. Note that,
if $D$ and $D'$ have non intersecting supports (which is easy to
realize), then $U_{D}$ and $U_{D'}$ cover $\C^{*}$ and $V_{D}$ and $V_{D'}$ 
cover $\Eq$.

\begin{lemma}
$F_{D,D'}(U)$ is a section of the sheaf $\La_{I}(M_{0})$ over $V_{D,D'}$.
\end{lemma}

\Pr 
One only has to prove that the upper-diagonal part of $F = F_{D,D'}(U)$
is flat. But its rectangular blocks satisfy:
$\l(\sq F_{i,j}\r) (z^{-\mu_{j}} A_{j}) = 
(z^{-\mu_{i}} A_{i}) F_{i,j}$.
This is a pure system of slope $\mu_{i,j} > 0$, hence $F_{i,j}$ 
is indeed flat by lemma \ref{lemma:flatness}.
\hfill $\Box$ \\

We now call $\U$ (resp. $\V$) the covering of $\C^{*}$ (resp. 
of $\Eq$) by the open subsets $U_{D}$ (resp. $V_{D}$), where 
$D$ runs among all the allowed summation divisors for $A_{0}$, 
$M_{0}$. The following is immediate.

\begin{cor}
The family $(F_{D,D'}(U))_{D,D'}$ is a Cech cocycle of the sheaf 
$\La_{I}(M_{0})$
for the covering $\V$ of $\Eq$.
\end{cor}

\hfill $\Box$

\begin{prop} \label{prop:cocycles-privilegies}
Mapping $A_{U}$ to the cocycle $(F_{D,D'}(U))$ defines a one-to-one
mapping:
$$
\F(M_{0}) \hookrightarrow Z^{1}(\V,\La_{I}(M_{0})).
$$
\end{prop}

\Pr
If $A_{U}$ is analytically equivalent to $A_{V}$, we have 
$A_{V} = F[A_{U}]$ for a unique $F \in \G(\Ka)$ and it is
clear (by unicity) that $F_{D}(V) = F F_{D}(U)$ for all
allowed divisors $D$, whence $F_{D,D'}(V) = F_{D,D'}(U)$
by immediate computation. This shows that the above mapping
is well defined. \\

Conversely, just assume that $F_{D,D'}(V) = F_{D,D'}(U)$
for two allowed divisors with non intersecting supports.
This equality gives $F_{D}(U,V) = F_{D'}(U,V)$. Hence, 
both sides are holomorphic over $U_{D} \cup U_{D'} = \C^{*}$
(near $0$), and we already saw in proposition \ref{prop:tout-dans-les-poles} 
that this implies the analytic equivalence of $A_{U}$ and $A_{V}$.
\hfill $\Box$ \\

We now come to the second fundamental result of this paper.

\begin{thm} \label{thm:q-Malgrange-Sibuya}
The above mapping yields a bijective correspondance:
$$
\F(M_{0}) \simeq H^{1}(\Eq,\La_{I}(M_{0})).
$$
\end{thm}

\Pr
Let $A_{U}$ and $A_{V}$ have the same image in
$H^{1}(\Eq,\La_{I}(M_{0})$, hence in $H^{1}(\V,\La_{I}(M_{0})$
(in Cech cohomology, the $H^{1}$ of a covering embeds into the
direct limit). There is, for each allowed divisor $D$, a matrix
$G(D) \in \G(\O(U_{D})) \cap Aut(M_{0})$ in such a way that:
$$
\forall D,D' \;:\;
F_{D,D'}(U) = \l(G(D)\r)^{-1} F_{D,D'}(V) G(D').
$$
One draws that $F_{D}(V) G(D) \l(F_{D}(U)\r)^{-1}$ does not depend
on $D$, so that it is holomorphic on $\bigcup U_{D} = \C^{*}$; call
$\Phi$ their common value. As a gauge transformation, it sends $A_{U}$
to $A_{V}$. By proposition \ref{prop:tout-dans-les-poles}, $A_{U}$
and $A_{V}$ are analytically equivalent, which proves the injectivity. \\

We now prove that our mapping from $\F(M_{0})$ to 
$H^{1}(\V,\La_{I}(M_{0}))$ is onto. For that, we 
take $(\Phi_{D,D'})_{D,D'} \in Z^{1}(\V,\La_{I}(M_{0}))$.
By definition of the sheaf $\La_{I}(M_{0})$, each component
$\Phi_{D,D'}$ is an element of $\G(\O(U_{D,D'}))$, so that
our cocycle can be considered as describing a vector bundle
over $\C^{*}$, trivialized by the covering $\U$ and with
structural group in $\G$. By \cite{Rohrl62}, theorem 1.0
(see also \cite{Rohrl57}) it is trivial in the following 
sense: there is, for each $D$, a $\Phi_{D} \in \G(\O(U_{D}))$
in such a way that, for all $D,D'$, 
$\Phi_{D,D'} = \Phi_{D}^{-1} \Phi_{D'}$. 
Since the $\Phi_{D,D'}$ are automorphisms of $A_{0}$,
the $\Phi_{D}[A_{0}]$ are all equal to a same matrix
$A_{U'}$. Moreover, this is holomorphic over $\C^{*}$.
By corollary \ref{cor:Rwg=Ka}, there is a 
$\Psi \in \G(\Rwg)$ such that $A_{U} = \Psi[A_{U'}]$
is meromorphic at $0$. Then $\Psi \Phi_{D}$ is holomorphic
on $(U_{D},0)$ and sends $A_{0}$ to $A_{U}$. Put
$G_{D} = F_{D}(U)^{-1} \Psi \Phi_{D}$. This is a section
of $\La_{I}(M_{0})$ over $V_{D}$. The equalities
$F_{D}(U) G_{D} = \Psi \Phi_{D}$ entail
$\Phi_{D}^{-1} \Phi_{D'} = G_{D}^{-1} F_{D,D'}(U) G_{D'}$,
that is, the cocycle $(\Phi_{D,D'})_{D,D'}$ is equivalent
to the cocycle $(F_{D,D'}(U))_{D,D'}$ which ends the proof
of our statement. \\

There remains to check that the natural mapping from
$H^{1}(\V,\La_{I}(M_{0}))$ to $H^{1}(\Eq,\La_{I}(M_{0}))$
is onto (we already said it was one-to-one). This is the 
content of proposition \ref{prop:bon-recouvrement}, to
be proved after the discussion on the $q$-Gevrey filtration 
of the Stokes sheaf the in next chapter.
\hfill $\Box$

%%%%%%%%%%%%%%%%%%%%%%%%%%%%%%%%%%%%%%%%%%%%%%%%%%%%%%%%%%%%%%%%%%%%%%%%%%%%%

% 4

\section{The $q$-Gevrey filtration on the Stokes sheaf}

The sources of inspiration for the contents of this chapter are
\cite{RamisGevrey}, \cite{BV}, \cite{MLR}, \cite{Varadarajan} and
\cite{Deligne86}. 

% 4.1

\subsection{The filtration for the Stokes sheaf of a tamely irregular module}

We stick to the conventions of section \ref{subsection:devissage},
in particular, the notations of \ref{equation:forme-canonique} and
\ref{equation:forme-canonique-diagonale}.

\subsubsection*{Conditions of flatness}

Let $F$ be a section of the sheaf $\La_{I}(M_{0})$. Then, 
for $1 \leq i < j \leq k$, the block $F_{i,j}$ is solution
of the equation:
$$ 
\sq F_{i,j} \l(z^{-\mu_{j}} A_{j}\r) = \l(z^{-\mu_{i}} A_{i}\r) F_{i,j}.
$$
>From this and lemma \ref{lemma:flatness}, we draw that $F_{i,j}$
is $(\mu_{i} - \mu_{j})$-flat and that, if it is $t$-flat for some
$t > \mu_{i} - \mu_{j}$, then it vanishes. \\

We now introduce a filtration of the Stokes sheaf and a filtration
of the sheaf of its Lie algebras. For real nonnegative $t$, we call 
$\la_{I}^{t}(M_{0})$ the subsheaf of $\la_{I}(M_{0})$ made of $t$-flat 
sections and $\La_{I}^{t}(M_{0})$ the subsheaf $I_{n} + \la_{I}^{t}(M_{0})$
of $\La_{I}(M_{0})$. The latter is a sheaf of unipotent subgroups,
while the former is the sheaf of its Lie algebras (see the discussion
in section \ref{subsection:irregular-theory}). Both filtrations
are decreasing and exhaustive (the $0$-term is the total sheaf, 
the $t$-term is the trivial sheaf for $t > \mu_{1} - \mu_{k}$). \\

>From the previous argument, we see that $\la_{I}^{t}(M_{0})$
has a very simple concrete description in terms of matrices: 
its sections have non trivial blocks only over the ``curved
over-diagonal'' consisting of those $(i,j)$-blocks such that
$\mu_{i} - \mu_{j} = t$. There is a similar description for
$\La_{I}^{t}(M_{0})$ (taking in account the block-diagonal
of identities). We shall however have use for a more intrinsic 
definition of these filtrations. We first describe the filtration 
of $\la_{I}(M_{0})$. 

\begin{prop} \label{prop:Lie=bundle} 
(i) The sheaf $\la_{I}(M_{0})$ is the sheaf of sections of the vector
bundle associated (see section \ref{subsection:devissage})
to the tamely irregular module $\underline{End}^{> 0}(M_{0})$. \\
(ii) The above filtration on the sheaf $\la_{I}(M_{0})$ is the
decreasing filtration associated to the graduation inherited from
$\underline{End}^{> 0}(M_{0})$.
\end{prop}

\Pr
We have $M_{0} = P_{1} \oplus \cdots \oplus P_{k}$, whence:
$$
\underline{End}(M_{0}) = 
\bigoplus_{1 \leq i,j \leq k} \underline{Hom}(P_{j},P_{i}).
$$
The internal Hom $\underline{Hom}(P_{j},P_{i})$ is a pure module of
slope $\mu_{i} - \mu_{j}$. Therefore, $\underline{End}^{> 0}(M_{0})$
is the sum of those $\underline{Hom}(P_{j},P_{i})$ such that
$\mu_{i} > \mu_{j}$, \emph{i.e.} $i < j$. \\

On the other hand, the vector bundle associated to 
$\underline{Hom}(P_{j},P_{i})$ has as sections on an open subset 
$V$ of $\Eq$ the morphisms from $P_{j}$ to $P_{i}$ that are holomorphic
on the preimage $U$ of $V$ in $\C^{*}$. This implies that the sheaf
of sections of the vector bundle associated to the module
$\underline{End}^{> 0}(M_{0})$ is indeed $\la_{I}(M_{0})$;
that it is the direct sum of its subsheaves $\la_{I}^{(t)}(M_{0})$, 
where  $\la_{I}^{(t)}(M_{0})$ is the sheaf of sections of the pure
vector bundle associated to the pure module
$$
\underline{End}^{(t)}(M_{0}) = 
\bigoplus_{\mu_{i} - \mu_{j} = t} \underline{Hom}(P_{j},P_{i});
$$
and that $\la_{I}^{t}(M_{0})$ is the direct sum of the
$\la_{I}^{(t')}(M_{0})$ for all $t' \geq t$.
\hfill $\Box$ \\

Actually the whole structure only depends on the filtrations
and the properties of internal Homs, so that it can be extended
to an arbitrary tannakian category.

\begin{prop}
Let $t$ be a nonnegative integer. \\
(i) $\La_{I}^{t}(M_{0})$ is a sheaf of normal subgroups 
of $\La_{I}(M_{0})$. \\
(ii) The map $f \mapsto 1 + f$ induces an isomorphism: 
$$
\lambda_{I}^{(t)}(M_{0}) \simeq 
\frac{\La_{I}^{t}(M_{0})}{\La_{I}^{t+1}(M_{0})}.
$$
\end{prop}

\Pr
Actually, these are purely algebraic properties: for a nilpotent 
two sided ideal $I$ of a non commutative algebra $A$, the subgroups
$1 + I^{t}$ of the unit group are normal and their successive
quotients are isomorphic to the quotient modules $I^{t}/I^{t+1}$.
\hfill $\Box$ \\

We now are in position to reconstruct the Stokes sheaf by successive
exact sequences:
\begin{equation} \label{equation:suite-exacte}
1 \rightarrow \La_{I}^{t+1}(M_{0}) \rightarrow
\La_{I}^{t}(M_{0}) \rightarrow \lambda_{I}^{(t)}(M_{0})
\rightarrow 0.
\end{equation}

Note also that, again from general algebraic considerations,
we have a sequence of central extensions:
\begin{equation} \label{equation:extension-centrale}
0 \rightarrow \lambda_{I}^{(t)}(M_{0}) \rightarrow
\frac{\La_{I}(M_{0})}{\La_{I}^{t+1}(M_{0})} \rightarrow
\frac{\La_{I}(M_{0})}{\La_{I}^{t}(M_{0})} \rightarrow 1.
\end{equation}

\subsection{Cohomological consequences}

\begin{lemma}
Let $V$ be a proper open subset of $\Eq$.
Then $H^{1}(V,\La_{I}(M_{0}))$ is trivial.
\end{lemma}

\Pr
We apply theorem I.2 of \cite{Frenkel} to the exact sequence
\ref{equation:suite-exacte}. This gives an exact sequence of
pointed sets:
$$
H^{1}(V,\La_{I}^{t+1}(M_{0})) \rightarrow 
H^{1}(V,\La_{I}^{t}(M_{0})) \rightarrow 
H^{1}(V,\la_{I}^{t}(M_{0})).
$$
If the extreme terms are trivial, so must be the central one
(this, by the very definition of an exact sequence of pointed
sets). The rightmost term is the first cohomology group of
a vector bundle (after proposition \ref{prop:Lie=bundle})
over an open Riemann surface. Such a bundle being a trivial
bundle, its $H^{1}$ is trivial. The leftmost term is trivial
for $t > \mu_{1} - \mu_{k}$. By descending induction, the
inner term is trivial for all $t$, hence for $t = 0$.
\hfill $\Box$

\begin{prop} \label{prop:bon-recouvrement}
The covering $\V$ is good.
\end{prop}

This means that the map from $H^{1}(\V,\La_{I}(M_{0}))$
to $H^{1}(\Eq,\La_{I}(M_{0}))$ is an isomorphism. After
\cite{BV}, cor. 1.2.4, p. 113, this follows from the
lemma.
\hfill $\Box$ \\

Note that this ends the proof of theorem \ref{thm:q-Malgrange-Sibuya}.

\subsection{The Stokes sheaf of a general module}

We briefly sketch here how the previous results extend
to the Stokes sheaf of a module $M$. We take $M$ in
the formal class of $M_{0}$ and identify it with
$(\Ka^{n},A_{U})$, according to the conventions of
section \ref{subsection:devissage}. \\

The mapping $\Phi \mapsto F_{D}(U) \Phi F_{D}(U)^{-1}$
defines an isomorphism from $\La_{I}(M_{0})$ to 
$\La_{I}(M)$ over $V_{D}$. Therefore, the two sheaves
are locally isomorphic. Actually, the latter is obtained
from the former by the operation of twisting by the
cocycle $(F_{D,D'}(U))_{D,D'}$, described in \cite{Frenkel},
prop. 4.2 (also see \cite{BV}, II.1 or \cite{Varadarajan}, 
pp. 30-31). According to the same references, their $H^{1}$
are isomorphic. Moreover, the same operations provide a local
isomorphism of the sheaves of Lie algebras, so that $\la_{I}(M)$
is a vector bundle. Last, these isomorphisms preserve the
filtrations by levels of flatness. \\

\begin{rmk} \label{rmk:vector-bundle}
One should also note that the mapping $X \mapsto F_{D}(U) X$
defines an isomorphism from the space of solutions of $A_{0}$
holomorphic over $(U_{D},0)$ to the same space for $A_{U}$.
This means that their sheaves of solutions are locally isomorphic.
Since the former is a vector bundle, so is the latter. This yields
an explicit way of associating a vector bundle with an arbitrary
module with integral slopes (for arbitrary slopes, see \cite{JSIRR},
\cite{JSISOMONO}).
\end{rmk}

%%%%%%%%%%%%%%%%%%%%%%%%%%%%%%%%%%%%%%%%%%%%%%%%%%%%%%%%%%%%%%%%%%%%%%%%%%%%%

% Biblio


\begin{thebibliography}{}

\bibitem{BV} \textsc{Babbitt D.G. and Varadarajan V.S.}
\emph{Local Moduli for Meromorphic Differential Equations},
Ast\'erisque, {\bf 169-170}, (1989).

\bibitem{Birkhoff1} \textsc{Birkhoff G.D.}
The generalized Riemann problem for linear differential equations and the
allied problems for linear difference and  $q$-difference
equations, \emph{Proc. Amer. Acad.}, {\bf 49}, (1913), 521-568.

\bibitem{Birkhoff3} \textsc{Birkhoff G.D. and Guenther P.E.}
Note on a Canonical Form for the Linear $q$-Difference System,
\emph{Proc. Nat. Acad. Sci.}, {\bf 27-4}, (1941), 218-222.

\bibitem{DF} \textsc{Deligne P.}
Cat\'egories Tannakiennes,
in \emph{Grothendieck Festschrift} (Cartier \& al. eds),
{\bf Vol. II},  Birkh\"{a}user, (1990).

\bibitem{DM} \textsc{Deligne P. and Milne J.}
\emph{Tannakian Categories},
in \emph{Hodge Cycles, Motives and Shimura Varieties} (Deligne \& al. eds),
\emph{Lecture Notes in Mathematics}, {\bf
$\text{n}^{\text{o}}$ 900}, (1989), Springer Verlag.

\bibitem{Deligne86} \textsc{Deligne P.}
Lettres \`a J.-P. Ramis de janvier et f\'evrier 1986.

\bibitem{DRSZ}
\textsc{Di Vizio L., Ramis J.-P., Sauloy J. and Zhang C.}
Equations aux $q$-diff\'erences,
\emph{Gazette des math\'ematiciens}, {\bf $\text{n}^{\text{o}}$
96}, (2003), 20-49.

\bibitem{Frenkel} \textsc{Frenkel J.}
Cohomologie non ab\'elienne et espaces fibr\'es,
\emph{Bull. S,M.F.}, {\bf 85}, (1957), 135-220.

\bibitem{MLR} \textsc{Loday-Richaud M.}
Stokes phenomenon, multisummability and differential Galois groups,
\emph{Ann. Inst. Fourier (Grenoble)}, {\bf 44-3},(1994), 849-906.

\bibitem{MZ} \textsc{Marotte F. and Zhang C.}
Multisommabilit\'e des s\'eries enti\`eres solutions formelles d'une
\'equation aux {$q$}-diff\'erences lin\'eaire analytique.
\emph{Ann. Inst. Fourier (Grenoble)}, {\bf 50-6},(2000),
1859-1890.

\bibitem{vdPS} \textsc{van der Put M. and Singer M.F}
\emph{Galois theory of difference equations},
\emph{Lecture Notes in Mathematics}, {\bf $\text{n}^{\text{o}}$
1666}, (1997), Springer Verlag.

\bibitem{RamisGevrey} \textsc{Ramis J.-P.}
D\'evissage Gevrey, in \emph{Journ\'ees singuli\`eres de Dijon},
Ast\'erisque {\bf 59-60}, (1978), 173-204.

\bibitem{RamisGrowth} \textsc{Ramis J.-P.}
About the growth of entire functions solutions to linear algebraic
$q$-difference equations,
\emph{Ann.  Fac.  Sciences Toulouse}, S\'erie 6, {\bf
I-1}, (1992), 53-94.

\bibitem{RamisJPRTraum} \textsc{Ramis J.-P.}
Fonctions $\theta$ et \'equations aux $q$-diff\'erences,
non publi\'e, (1990), Strasbourg.

\bibitem{RSZ} \textsc{Ramis J.-P., Sauloy J. and Zhang C.}
Local analytic classification of  irregular $q$-difference
equations, \emph{Article in preparation}.

\bibitem{RSZcras1} \textsc{Ramis J.-P., Sauloy J. and Zhang C.}
La vari\'et\'e des classes analytiques d'\'equations aux
$q$-diff\'erences dans une classe formelle.
C. R. Math., {\bf 338, 4}, (2004), 277-280.

\bibitem{RZ} \textsc{Ramis J.-P. and Zhang C.}
D\'eveloppement asymptotique $q$-Gevrey et
fonction th\^eta de Jacobi,
\emph{C. R. Acad. Sci. Paris}, {\bf Ser. I 335}, (2002), 899-902.

\bibitem{Rohrl57} \textsc{R\"{o}hrl H.}
Das Riemann-Hilbertsche Problem der Theorie der linearen
Differentialgleichungen,
\emph{Math. Ann.}, {\bf 133}, (1957), 1-25.

\bibitem{Rohrl62} \textsc{R\"{o}hrl H.}
Holomorphic fiber bundles over Riemann surfaces,
\emph{Bull. A.M.S.}, {\bf 68}, (1962), 125-160.

\bibitem{JSAIF} \textsc{Sauloy J.}
Syst\`emes aux $q$-diff\'erences singuliers r\'eguliers :
classification, matrice de connexion et monodromie,
\emph{Ann. Inst. Fourier (Grenoble)},
{\bf 50-4}, (2000), 1021-1071.

\bibitem{JSGAL} \textsc{Sauloy J.}
Galois theory of fuchsian $q$-difference equations,
\emph{Ann. Sci. \'Ecole Norm. Sup. (4)}, {\bf 36}, (2003), no 6, 925-968.

\bibitem{JSFIL} \textsc{Sauloy J., 2002.}
La filtration canonique par les pentes d'un module aux $q$-diff\'erences
et le gradu\'e associ\'e,
\emph{Ann. Inst. Fourier}, {\bf 54}, 1, (2004), 181-210.

\bibitem{JSIRR} \textsc{Sauloy J.}
Local Galois theory of irregular $q$-difference equations,
\emph{Article in preparation}.

\bibitem{JSISOMONO} \textsc{Sauloy J.}
Isomonodromy for $q$-difference equations,
\emph{Article to be submitted}.

\bibitem{Varadarajan} \textsc{Varadarajan V. S.}
Linear meromorphic differential equations: a modern point of view,
\emph{Bull. A.M.S.} {\bf 33-1}, (1996), 1-42.

\bibitem{Zha99} \textsc{Zhang C.}
D\'eveloppements asymptotiques {$q$}-{G}evrey et s\'eries
{$Gq$}-sommables.
\emph{Ann. Inst. Fourier (Grenoble)}, {\bf 49-1}:vi--vii, x,
(1999) 227--261.

\bibitem{Zha02} \textsc{Zhang C., 2002.}
Une sommation discr\`ete pour des \'equations
aux $q$-diff\'erences lin\'eaires et \`a coefficients
analytiques: th\'eorie g\'en\'erale et exemples,
in \emph{Differential Equations and the Stokes Phenomenon},(2002),
  B.L.J. Braaksma, G. Immink, M. van der Put and J. Top, editors,
World Scientific.

\end{thebibliography}
\end{document}